\documentclass[20p]{amsart}

\usepackage{graphics}
\usepackage{shadow}
\usepackage[latin1]{inputenc}
\usepackage[frenchb]{babel}

\theoremstyle{definition}

\theoremstyle{remark}

\numberwithin{equation}{section}

\newcommand{\be}{\begin{equation}}
\newcommand{\ee}{\end{equation}}

\newcommand{\bbR}{{\mathbb R}}

\newcommand{\p}{Poincaré\ }

\begin{document}

\title{Poincaré et les quanta}
\author[T.Paul]{ T. Paul\\ C.N.R.S. et D.M.A., Ecole Normale Supérieure}

\setcounter{page}{0}

\maketitle
\tableofcontents
\LARGE
\section{Introduction}
\p s'intéressa très tard à la théorie des quanta. Invité au Congrès
Solvay de 1911, il découvrit la problématique quantique et écrivit alors coup sur coup plusieurs
articles sur le sujet. Certains de ``vulgarisation" \cite{PO}, et deux articles
scientifiques \cite{PO1,PO2} dans lesquels il va s'efforcer de montrer que l'hypothèse de Planck ne peut
être contournée et remplacée par une théorie ``continue". Il  s'agît pour lui de se demander si la formule de Planck donnant la densité d'énergie du corps noir, et
dérivée depuis l'hypothèse des quanta (qui consiste, pour dire vite, à remplacer dans
un calcul une intégrale par une somme) peut être (ou non) déduite d'un modèle continu d'interaction entre rayonnement et molécules. Le résultat que \p lui-même
présentera comme négatif validera la théorie de Planck.

Ces deux articles sont peu cités de nos jours, mais ont  eu un retentissement important l'époque. En est la preuve le fait que deux articles dans le Volume 38
d'Acta Mathematica en hommage à \p y sont consacrés. De plus l'argument de \p sera repris et d'une certaine manière rendu plus
rigoureux dans le célèbre traité de Mécanique Statistique de Fowler en 1929 \cite{FO}. Leur (re)lecture après plus d'un siècle de Mécanique Quantique est intéressante
à plusieurs égards. Tout d'abord ils montrent, à travers une discussion parfaitement technique, l'affrontement entre discret et continu, en ce début de XXième siècle
 où la
problématique fondationnelle des mathématiques fait rage. Ensuite, et justement autour de cette dualité discret/continu, ils présentent de merveilleuses mathématiques.
On est pressent en particulier la masse de Dirac, et le lemme de la phase stationnaire, qui deviendra si important par la suite pour la Mécanique Quantique y est utilisé. 
Enfin et surtout ils exhibent clairement une part méthodologique de son héritage scientifique : l'intérêt pour les
résultats négatifs, pour les erreurs. depuis l'erreur dans le système des 3 corps jusqu'à son ``échec"  dans la tentative de remplacer l'hypothèse des quanta \p
n'aura cessé de décortiquer le négatif pour y faire naître le constructif.

Rappelons pour terminer cette introduction que la principale contribution de \p à la Mécanique Quantique nous semble consister en le deuxième Tome des ``Méthodes nouvelles de la Mécanique
Céleste". C'est avec les calculs de perturbations de la mécanique céleste que Max Born essaiera dans les années 20 (en le citant explicitement, voir \cite{mb}) de quantifier perturbativement l'hélium, avant dy plonger son jeune assistant W.
Heisenberg qui fera sortir de l'algèbre de convolution des séries de Fourier l'algèbre non-commutative des matrices (voir \cite{tp}).

\section{Chronologie}

\[
\begin{array}{cccc}
&Planck &&Poincar\acute{e}\\
&&&\\
1854& && Naissance\\
&&&\\
1858&Naissance&&\\
&&&\\
1889& &&Roi\ de\ Su\grave{e}de\\
&&&\\
1892& && M\acute{e}thodes\ nouvelles\\
&&&\\
1900&quanta&&\\
&&&\\
1905 &Photon& &\\
&&\\
1911 &Congr\grave{e}s&&Solvay\\
&&\\
1912 & && Sur\ la\ th\acute{e}orie\ des\ quanta\\
&&\\
1912 & && Mort\\
&&\\
1915 & Hommages&\grave{a} & Poincar\acute{e}\\
&&\\
1925/6& M\acute{e}canique\ quantique&&\\
&&\\
1929& & &Trait\acute{e}\ de\ Fowler\\
&&\\
1947&Mort&&\\

\end{array}
\]
\newpage
\section{La théorie des quanta en 1911/12}
1911 se situe dans une période un peu creuse pour la théorie quantique. L'hypothèse des quanta, énoncée par Planck en 1900 et généralisée par Einstein avec l'effet
photoélectrique en 1905, est certes bien acceptée, mais une véritable  dynamique quantique (Heisenberg 1925) et l'application \`a la théorie atomique (Bohr 1913) se
font attendre. De plus  Planck, qui ne sera jamais un supporter fanatique de la mécanique  quantique, est réticent (il semble qu'il soit toujours resté insatisfait de sa
loi qu'il aurait voulu dériver de principes strictement thermodynamiques), et l'on peut bien imaginer qu'Einstein est
fort occupé à inventer la relativité générale. Est-ce pour cette raison que Poincaré ne s'est jamais, auparavant, vraiment intéressé aux quanta? La réponse nous
semble justement se trouver dans l'article que nous allons décrire et commenter : le changement paradigmatique des quanta est grand, déjà à l'époque, et bien avant les
``problèmes conceptuels" de la mécanique quantique.

Car il s'agit de l'opposition continu/discret, et de l'idée d'une théorie physique qui puiserait ses sources hors du cadre ``naturel" des équations différentielles de la
culture classique. L'hypothèse des quanta consiste, dans notre langage moderne, à discrétiser une intégrale (ce que fait notre ordinateur), mais sans faire tendre, a
posteriori,  le
pas de discrétisation vers zéro.

\p est invité au Congrès Solvay de 1911 et il prend à bras le corps la problématique quantique: ne serait-il pas possible de dériver la même loi que Planck sans
l'hypothèse discrète?. \p y consacrera une note et un article,  et donnera une réponse négative. Il mourra en 1912, en 1913 Niels Bohr appliquera avec le succès que l'on sait la
théorie au modèle planétaire de l'atome, et en 1926 Erwin Schr\"odinger définira une équation aux dérivées partielles (donc plongée dans le {\it continu}\rm), qui n'aura
de solution raisonnable que pour des valeurs {\it discrètes}\rm\ des paramètres d'énergie : la boucle sera donc bouclée 15 ans après l'article de Poincaré que nous
allons présenter maintenant.

\section{L'article de Planck de 1900}
Nous n'allons pas ici rappeler la théorie du corps noir ni présenter en détail la contribution de Planck. Nous verrons dans la section suivante la 
dérivation ``critique" de la loi de Planck
par Poincaré lui-même. Disons seulement que le problème consiste à  trouver une formule d'interpolation pour la densité
d'énergie $u$ d'un corps ``noir", exprimée par rapport à la fréquence du rayonnement $\nu$, telle que \[
u(\nu)\sim \nu^2,\ \nu\to 0\]
et \[
u(\nu)\sim e^{-cte \nu},\  \nu\to\infty.
\]
Une formule simple est celle imaginée par Planck :
\[
u(\nu)=\frac{8\pi\nu^2}{c^3}\frac{h\nu}{e^{h\nu/kT}-1}.
\]
C'est cette loi que Planck  dérive en 1900 sous l'hypothèse des quanta. Plus
exactement il  montre que, si l'on suppose que l'énergie de chaque oscillateur
est un multiple entier d'une certaine quantité $\epsilon$, donc de la forme $P\epsilon$, $P$ entier, un calcul combinatoire 
donne que l'entropie $S$ de chaque oscillateur est une fonction (explicite) de
$\frac U \epsilon$, où $U$ est l'énergie de chaque oscillateur. Plus précisément 
\be\label{entr}
S=k\left(\left(1+\frac U \epsilon\right)log\left(1+\frac U \epsilon\right)-\frac U \epsilon log\frac U \epsilon\right)
\ee
Utilisant la loi de Wien concernant la température 
\[
\frac 1 {kT}=\frac{dS}{dU},
\]
où $k$ est la constante de Boltzmann
 et un peu d'analyse dimensionnelle qui donne que  
$\epsilon =h.\nu$, où $h$ a la dimension d'une action et $\nu$ l'inverse d'un temps,
on peut  intégrer (\ref{entr}). On obtient :
\[
u:=\frac{8\pi\nu^2}{c^3}U=\frac{8\pi\nu^2}{c^3}\frac{h\nu}{e^{h\nu/kT}-1}.
\]
C'est la loi de Planck.
\section{Les articles de Poincaré}
Poincaré va reprendre le calcul en partant d'un densité d'énergie $w$ et montrer
tout d'abord que si $w$ est ce que l'on appellerait aujourd'hui un ``peigne de
Dirac" sur les entiers ($\times h$) il retrouve la loi de Planck. C'est donc en
quelque sorte une vision plus mathématique , ou plutôt moins thermodynamique, que
Planck ; mais il va surtout faire le chemin inverse, et ``montrer" que la forme de
la fonction $u$  détermine $w$ de façon univoque  et donc que seule
l'hypothèse des quanta donne le bon résultat.

Poincaré commence par montrer que la densité d'énergie 
$$W(\eta_1,\dots,\eta_n;\xi_1,\dots,\xi_p)$$ (on écrit à l'époque 
$W(\eta_1,\dots,\eta_n;\xi_1,\dots,\xi_p)d\eta_1...d\eta_nd\xi_1...d\xi_p$ de n oscillateurs (identiques) d'énergies
$\eta_1\dots\eta_n$ et p molécules 
(identiques) d'énergie $\xi_1\dots\xi_p$ 
peut prendre  la forme :
\[
W(\eta_1,\dots,\eta_n;\xi_1,\dots,\xi_p)=\Pi_{i=1}^{n}w(\eta_i)=w(\eta_1)...w(\eta_n)
\]
pour une certaine fonction $w(\eta)$. C'est cette fonction $w$ qu'il va montrer devoir avoir une primitive discontinue.

Pour cela il considère tout d'abord, dans $\bbR^{n+p}$ la surface d'énergie $S_h=\{(\eta_1,\dots,\eta_n;\xi_1,\dots,\xi_p)/ \eta_1+\dots+\eta_n+\xi_1+\dots+\xi_p=h\}$.
Puis les 3 intégrales :
\[
\begin{array}{lcl}
I&=&\int_S w(\eta_1)...w(\eta_n)d\eta_1...d\eta_nd\xi_1...d\xi_p \\&&\\
I'&=&\int_S(\eta_1+\dots+\eta_n) w(\eta_1)...w(\eta_n)d\eta_1...d\eta_nd\xi_1...d\xi_p \\&&\\
I"&=&\int_S(\xi_1+\dots+\xi_p)w(\eta_1)...w(\eta_n)d\eta_1...d\eta_nd\xi_1...d\xi_p\end{array}
\]
On définit alors les énergies moyennes des résonnateurs et des molécules, $X$ et $Y$, par  :
\[nYI=I'\ \mbox{et}\ pXI=I''\]
Si l'on pose maintenant 
\be\label{phi}
\int_{\eta_1+\dots+\eta_n=x}w(\eta_1)...w(\eta_n)d\eta_1...d\eta_n=\phi(x),
\ee
un calcul simple donne
\[
I=\frac 1 {(p-1)!}\int_0^h(h-x)^{p-1}\phi(x)dx
\] De la même façon
\[
I'=\frac 1 {(p-1)!}\int_0^hx(h-x)^{p-1}\phi(x)dx\]et \[ I''=\frac 1 {(p-1)!}\int_0^h(h-x)(h-x)^{p-1}\phi(x)dx,
\]
et finalement :
\[
nY=\frac{\int_0^hx(h-x)^{p-1}\phi(x)dx}{\int_0^h(h-x)^{p-1}\phi(x)dx}\ \mbox{et}\ pX=\frac{\int_0^h(h-x)^p\phi(x)dx}{\int_0^h(h-x)^{p-1}\phi(x)dx}.
\]

Poincaré montre alors que l'on obtient la formule de Planck si l'on choisit pour $w$ la ``fonction" définie, à partir de $\epsilon>0$, de la façon suivante :

\[
w(\eta)=0 \ \ \ \mbox{si}\ \ \ k\epsilon<\eta<k\epsilon+\mu\ (k\in\mathbb N)\ \forall\mu,\ \ 0<\mu<\epsilon
\]
et
\[
\int_{k\epsilon}^{k\epsilon+\mu}w(\eta)d\eta=1,\ \forall\mu,\ \ 0<\mu<\epsilon.
\]

Mais Poincaré de s'arrête pas là. Comme il l'a dit dans l'introduction de son article \cite{PO2}, il veut absolument s'assurer que ce choix de $w$ 
est le seul qui donne la loi de Planck, et surtout qu'aucun choix ``continue" ne redonne le même résultat. Nous allons pr\'esenter en m\^eme temps la d\'erivation directe et la
solution du probl\`eme inverse.

\p va montrer que le rapport $\frac Y X$ 
détermine $w$.
Puisqu'il a montré que le choix discontinu donne la bonne formule il aura gagné (mais dira qu'il aura perdu).

L'argument  utilise  l'analyse complexe et nous allons le décrire rapidement. On introduit tout d'abord la transformée de Laplace (que \p appelle ``intégrale de Fourier") de $w$  :
\[
\Phi(\alpha)=\int_0^\infty w(\eta)e^{-\alpha\eta}d\eta.
\]
 $w$ (et donc $\phi$) est déterminé par $\Phi$ grâce à la formule d'inversion de la transformé de Laplace :
 \[
 w(\eta)=\frac 1 {2\pi i}\int_L\Phi(\alpha)e^{\alpha\eta}d\alpha,
 \]
où $L$ est une droite complexe parallèle à l'axe imaginaire dans le demi-plan $\{\Re\alpha>0\}$.

D'autre part, par définition (\ref{phi}) de $\phi$, on a :
\[
\left(\Phi(\alpha)\right)^n=\int_0^\infty \phi(x)e^{-\alpha x}dx,
\]
puisque
\begin{eqnarray}
\left(\Phi(\alpha)\right)^n&=&\left(\int_0^\infty w(\eta)e^{-\alpha\eta}d\eta\right)^n\nonumber\\
&=&\int w(\eta_1)\dots w(\eta_n)e^{-\alpha(\eta_1+\dots+\eta_n)}d\eta_1\dots d\eta_n\nonumber\\
&=&\int_{x=0}^\infty\left(\int_{\eta_1+\dots+\eta_n=x}w(\eta_1)...w(\eta_n)d\eta_1...d\eta_n\right)e^{-\alpha x}dx\nonumber\\
&=&\int_0^\infty \phi(x)e^{-\alpha x}dx.\nonumber
\end{eqnarray}
D'où :
\be\nonumber
\phi(x)=\frac 1 {2\pi i}\int_L\left(\Phi(\alpha)\right)^ne^{\alpha\eta}d\alpha,
\ee
et donc, après la renormalisation $x=n\omega,\ h=n\beta$, et avoir posé $$\Theta(\alpha,\omega):=\Phi(\alpha)e^{\alpha\omega}(\beta-\omega)^{\frac p n},$$ on obtient :
\[
nY=\frac{\frac{n^{p+1}}{2\pi i}\int_0^\beta\int_L\frac{\omega}{\beta-\omega}\Theta^nd\omega d\alpha}{\int_0^h(h-x)^{p-1}\phi(x)dx}
\]
et
\be\label{mnb}
pX=\frac{\frac{n^{p+1}}{2\pi i}\int_0^\beta\int_L\Theta^nd\omega d\alpha}{\int_0^h(h-x)^{p-1}\phi(x)dx}\ .
\ee

Vient ensuite l'utilisation du fait que le système physique con-\\tient un grand nombre d'oscillateurs (on a un champ de rayonnement). Poincaré fait tout d'abord 
 l'hypothèse que $\Theta$ atteint un maximum à $\omega=\omega_0$ et $\alpha=\alpha_0$. 

Puisque $n$ est très grand et que les intégrales précédentes
font intervenir $\Theta^n$, il en déduit, en choisissant $L=\{\Re z=\alpha_0\}$ que :
\be\label{frac}
\frac{nY}{pX}\sim \frac{\omega_0}{\beta-\omega_0}.
\ee
 mais puisque, par ailleurs, $nY+pX=h=n\beta$, on déduit que 
 \[
 Y=\omega_0\ \ \ \mbox{et}\ \ \ X=\frac{\beta-\omega_0}k.
 \]
 En fait Poincaré utilise le théorème de la phase stationnaire qui s'applique ici, puisque, à une constante près, les intégrales précédentes font intervenir la
 quantité $\Theta(\alpha,\omega)^n=e^{n\mbox{log}(\Theta(\alpha,\omega)}$. En effet le numérateur de (\ref{mnb}), par exemple,s'écrit :
 \[ 
\frac{n^{p+1}}{2\pi i}\int_0^\beta\int_Le^{n\mbox{log}(\Theta(\alpha,\omega)}d\omega d\alpha.
\]
 Le théor\`eme de la phase stationnaire dit alors que l'intégrale précédente se ramène, dans la limite $n\to\infty$, à la contribution de l'intégrant aux 
 ``points critiques" $\alpha_0$ et $\omega_0$  résolvant les équations :
 \[
 \frac{\partial\mbox{log}(\Theta(\alpha,\omega)}{\partial\alpha}=\frac{\partial\mbox{log}(\Theta(\alpha,\omega)}{\partial\omega}=0.
 \]
 ce qui donne 
 \[
 \frac{\Phi'(\alpha_0)}{\Phi(\alpha_0)}+\omega_0=0
 \]
 et
 \[
 \alpha_0-\frac k{\beta-\omega_o}=0,
 \]
 d'où (\ref{frac}).

 On voit donc que $X=\frac 1 \alpha_0$. 
 
 Utilisant le fait (Boltzmann) que l'énergie moyenne d'une molécule est proportionelle à la température (absolue) $T$ on en déduit que :
 \[
 \alpha_0=\frac C T=\frac 1 {kT}\]
 où $C$ est une constante, et donc $Y=\omega_0=- \frac{\Phi'(\alpha_0}{\Phi(\alpha_0)}$ nous donne  l'énergie moyenne d'un r\'esonnateur en fonction de la température (qui
 d'ailleurs sera indépendant du rapport $\frac n p$).  
  Il suffit alors de remarquer que, si la fonction $w$ satisfait l'hypothèse de Planck, alors :
  \begin{eqnarray}
  \Phi(\alpha)&=&\int\sum_k \delta(\eta-k\epsilon)e^{-\alpha\eta}d\eta\nonumber\\
  &=&\sum_k e^{-k\epsilon\alpha}\nonumber\\
  &=&\frac{e^{-\epsilon\alpha}}{1-e^{-\epsilon\alpha}}\nonumber\\
  &=&\frac 1 {e^{\epsilon\alpha}-1}.\nonumber
  \end{eqnarray}
  Et donc que, pour la densité $U$ de la section précédente,
  \[
  U=Y=- \frac{\Phi'(\alpha_0}{\Phi(\alpha_0)}=\frac{\epsilon}{e^{\frac\epsilon{kT}}-1}=\frac{h\nu}{e^{\frac{h\nu}{kT}}-1},
  \]
  ce qui est la loi de Planck.

 Mais supposons maintenant que l'on connaisse, {\it pour toutes les températures}\rm, l'énergie moyenne d'un r\'esonnateur. Alors on en déduira les valeurs de la fonction
 $\frac{d\mbox{log}\Phi(\alpha)}{d\alpha}$ pour tous les $\alpha>0$. Soit, à une constante multiplicative près, les valeurs de $\Phi$ sur l'axe réel, et par continuation analytique (nous
 dit Poincaré) à tout le demi-plan. Et donc, par la formule d'inversion précédemment énoncée, la fonction $w$ sera déterminée (toujours à une constante multiplicative près).
  
 Poincaré déduit alors de ce résultat d'unicité la nécessitée de l'hypothèse des quanta.
 
 Finalement Poincaré montre encore que, sans se servir de la loi  de Planck, et en confrontant sa formule donnant l'énergie d'un r\'esonnateur à celle du rayonnement du corps
 noir, $w$ doit vérifier, sous peine que l'énergie totale diverge, que  :
 \[
 \lim_{\eta_0\to 0}\int_0^{\eta_0}w(\eta)d\eta\neq 0.
 \]
  C'est à dire que $w$ doit être singulière.

\section{La masse de Dirac}

La condition que $w$ soit "choisie" régale à 0 sauf pour
$k\epsilon<\eta<k\epsilon+\delta$", où "$\delta$ est un infiniment petit" et telle que :
\[
\int_{k\epsilon}^{k\epsilon+\delta}w(\eta)d\eta=1
\]
s'exprime en langage moderne en disant que $w$ est une distribution, somme de masses de
Dirac sur le réseau des $\{k\epsilon\ k=1,\ 2\ \dots\}$ (pour être rigoureux, la condition intégrale doit être
$\int_{k\epsilon-\delta}^{k\epsilon+\delta}w(\eta)d\eta=1$).

Ceci n'est jusque là qu'une red\'erivation, mais sous forme intégrale et non
plus combinatoire, de la formule de Planck.


\section{Que nous apprend presque un siècle après l'article de Poincaré?}
Une nouveauté contenue dans cet article tient justement à écrire une somme comme une intégrale ... d'une fonction singulière. Les mathématiques du
XXième siècle nous ont appris à faire cela. Les distributions de Laurent Schwartz en sont l'outil parfait : par intégration (cette possibilité tellement ``continue") 
elles nous font apparaître des points (entités de mesure nulle). Jean Leray, dans son magnifique résultat d'existence des équations de Navier-Stokes \cite{JL} 
définira lui
aussi, en 1934, des solutions faibles pour les équations de l'hydrodynamique censées représenter des flots de particules {\it ponctuelles}\rm \ en mouvement. Mais l'idée de
chercher un intégrant singulier qui donnerait, par intégration, une somme discrète est tout à fait étonnante en 1911.

Ensuite l'article de Poincaré nous montre un exemple merveilleux de confrontation entre mathématiques et physique. L'hypothèse des quanta est considérée comme physique,
comme un changement paradigmatique dans la physique. Elle n'a pas à être justifiée en dehors de son statut d'hypothèse d'un ``élément de réalité" comme le dit 
lui-même \p. Mais
\p lui donne un autre statut: celui d'une nécessité mathématique, dans un certain cadre conceptuel, en vue de la dérivation d'une formule donnée. C'est
cet argument-là qui me semble expliquer que deux des articles (signés par Lorentz et de Planck!) du Volume d'Acta Mathematica  de 1915 en hommage à Poincaré 
soient consacrés à sa contribution aux quanta \cite{LO, PLA}, et que quelques 27 ans plus tard Fowler y consacre quelques pages dans son traité de Mécanique Statistique
(dans lesquelles il va rendre ``rigoureux" le raisonnement de \p).

Enfin cet article nous montre l'attitude de \p face aux résultats négatifs. car au fond il nous démontre un résultat positif : l'hypothèse des quanta est nécessaire à
la dérivation de la formule de Planck. Mais il nous présente son raisonnement par la négative, comme un échec (``je suis arrivé à un résultat négatif"
\cite{PO1}, ``je dois dire tout de suite que j'ai été conduit à répondre négativement à la question posée... "\cite{PO2}). 

Cet exemple, à la fin de sa vie, vient corroborer celui des 3 corps (``j'attire l'attention sur les résultats négatifs présents dans ce mémoire" \cite{POI}) du
début et nous montre \p dans une véritable activité positive de déconstruction du négatif.

\end{document}